\documentclass{amsart}
 \usepackage{amsmath}
 \usepackage{epsfig} 
 \usepackage{graphics} 

\newtheorem{theorem}{Theorem}
\newtheorem{corollary}{Corollary}
\newtheorem*{main}{Main Theorem}
\newtheorem{lemma}{Lemma}

\newtheorem{proposition}{Proposition}

\newtheorem{problem 1}{Problem 1}
\newtheorem{problem 2}{Problem 2}
\newtheorem{problem 3}{Problem 3}
\newtheorem{definition}{Definition}

\newcommand{\bea}{\begin{eqnarray*}}
\newcommand{\eea}{\end{eqnarray*}}

\title[LAMINATED CURRENTS]
      {LAMINATED CURRENTS}

\begin{document}

\begin{abstract}
In this paper we prove  the equivalence of two definitions
of laminated currents.
 \end{abstract}

\maketitle

\centerline{\scshape  John Erik Forn\ae ss\footnote{The first author is supported by an NSF grant. Keywords: Approximation, currents, test functions.
2000 AMS classification. Primary: 57R30, Secondary: 32U40}, Yinxia Wang and Erlend Forn\ae ss Wold}

\


\medskip

 \medskip


\section{Introduction}

 Let $K$ be a relatively closed subset of the bidisc $\Delta^2(z,w)=\{(z,w); |z|,|w|<1\}.$
 We suppose that $K$ is a disjoint union of holomorphic graphs,
 $w=f_\alpha(z),$ where $f_\alpha$ is a holomorphic function on the unit disc
 with $f_\alpha(0)=\alpha$ and $|f_\alpha(z)|<1$. We let $\mathcal{L}$ denote the lamination of  $K$. 
\

There are two notions of laminated currents
that we will discuss, \cite{FS2005}. Let $T$ be a positive closed $(1,1)$ current supported on $K.$
We assume that $T$ is the restriction of a positive closed current defined on a neighborhood of
 $\overline{\Delta}^2.$ We denote by
 $[V_\alpha]$ the current of integration along the graph of $f_\alpha.$
 Let $\lambda$ denote a continuous $(1,0)$ form which at $(z,f_\alpha(z))$ equals a non zero
 multiple of
 $dw-f'_\alpha(z)dz.$

 \begin{definition}
 We say that $T$ is weakly directed by the lamination ${\mathcal L}$  if
 $\lambda \wedge T=0$ for any such $\lambda.$
 \end{definition}

 \begin{definition}
 We say that $T$ is directed by ${\mathcal L}$ if there is a positive measure $\mu$ so that
 $T=\int_\alpha  [V_\alpha]d\mu(\alpha).$
 \end{definition}

Our main result is
\begin{main}
The current $T$ is directed if and only if it is weakly directed.
\end{main}

We note that this is a result by Sullivan in the case of the lamination being smooth,
i.e. the graphs vary smoothly with $\alpha,$ \cite{S1976}.  The part of Sullivan's proof that does not go through automatically in the non smooth case is a certain approximation step, and so in the present article we are concerned with approximation of partially smooth functions.  In \cite{FWW2007} the authors proved such an approximation theorem in the case of laminations in $\mathbb R^2$ and in $\mathbb R^3.$  

In the last section we show that the main theorem breaks down for Riemann surface laminations in higher dimension.
\medskip

\section{Preliminary estimates for slopes of holomorphic graphs}

The case we are studying is given by a holomorphic motion, see \cite{AM2001} for an exposition
and further references. In this paper we need a basic estimate on slopes of the graphs.
For the benefit of the reader we include the details of this well known fact.

We start with a Lemma. Let $\Delta:=\{z\in \mathbb C; |z|<1\}$ denote the unit disc in $\mathbb C.$
We denote by $\mathcal O(\Omega)$ the space of holomorphic functions on $\Omega.$
Let $\|\cdot\|_\infty$ denote the sup norm. Set
$$
H^\infty=H^\infty (\Delta)= \{f\in \mathcal O(\Delta); \|f\|_\infty <\infty\}.
$$
 Also, if $0<C<\infty$ we set

 $$
H^\infty_C=H^\infty_C (\Delta)= \{f\in \mathcal O(\Delta); \|f\|_\infty <C\}.
$$

\begin{lemma} If
 $f\in H^\infty_1(\Delta) $ and $f(z)\neq 0\; \forall z \in \Delta,$
then $|f'(0)| \leq 2 |f(0)| \log \frac{1}{|f(0)|}.$
\end{lemma}

\begin{proof}
Pick a holomorphic function $f(z)$ on the unit disc such that $0 \neq |f(z)|<1$ for all $z \in \Delta.$
We can replace $f(z)$ by $e^{i\theta}f(z)$ for any real $\theta. $ This does not change
$|f(0)|$ and $|f'(0)|.$ Hence we can assume that $f(0)>0.$ 

We set $h(z):=\log f(z).$ Then $h(z)$ is a holomorphic function on the unit disc
and ${\mbox{Re}}(h(z))<0.$ We can also choose a branch of the logarithm so that
$\log (f(0))=-a<0.$ If $k(z)=\frac{h(z)}{a},$ then $k(z)$ is a holomorphic function on the
unit disc and $k(0)=-1, {\mbox{Re}}(k(z))<0.$ We define $L(w)=\frac{w+1}{w-1}.$ Then
$L(-1)=0$ and if ${\mbox{Re}}(w)<0$ then $|L(w)|<1.$
Then $\Gamma(z):=L(k(z))$ is a holomorphic function from the unit
disc to the unit disc. Moreover $\Gamma(0)=L(k(0))=L(-1)=0.$
Since $\Gamma(0)=0$
and $|\Gamma(z)|<1$ we can apply the Schwarz' Lemma. So we can conclude that
$|\Gamma'(0)| \leq 1.$
By the chain rule, $\Gamma'(0)=L'(k(0))k'(0)=L'(-1) k'(0).$
Since $L'(w)
=\frac{-2}{(w-1)^2}$ we get $\Gamma'(0)=\frac{-2}{(-1-1)^2} k'(0)$ and
therefore $k'(0)=-2 \Gamma'(0).$ Hence we get $|k'(0)| \leq 2.$
Since $k(z)=\frac{h(z)}{a},$ we next can conclude that $|k'(0)|=|h'(0)|/a.$
Hence $|h'(0)|=a|k'(0)|\leq a \cdot 2$ so $|h'(0)| \leq 2a.$
Next recall that $h(z)=\log f(z)$ so $f(z)=e^{h(z)}.$ Hence $f'(z)=e^{h(z)}h'(z).$
Therefore $f'(0)=e^{h(0)}h'(0)=f(0) h'(0).$ Hence $|f'(0)| \leq |f(0)||h'(0)|.$
This implies that $|f'(0)|\leq 2a|f(0)|. $ Now recall that $\log f(0)=-a.$
But we have set this up so that $\log f(0)=\log |f(0)|+i \arg f(0) $ is real valued. So
we have that $\log |f(0)|=-a$ i.e. $\log \frac{1}{|f(0)|}=a.$ Therefore
$|f'(0)| \leq 2a|f(0)| =2|f(0)| \log \frac{1}{|f(0)|}.$ This concludes the proof of the Lemma.

\end{proof}

\begin{corollary}\label{estimate}
Suppose that we have two functions $f$ and $g$ holomorphic on the
unit disk with  $f-g\in H^\infty_{1}(\Delta).$ Suppose that $f(z)
\neq g(z)$ for each $z \in \Delta.$ We then have the estimate
$|f'(z)-g'(z)| \leq 4 |f(z)-g(z)| \log \frac{1}{|f(z)-g(z)|}$ for
all $z\in \Delta, |z|<1/2.$
\end{corollary}

\begin{proof}
Pick $z, |z|<1/2.$ We define $G(w)= f(z+w/2)-g(z+w/2)$. Then $G(w)$
satisfies the conditions of Lemma 1. Hence $ |G'(0)|\leq
2{|G(0)|}\log{\frac{1}{|G(0)|}}$. Therefore,

$$
\frac{1}{2}|f'(z)-g'(z)|  \leq
2{|f(z)-g(z)|}\log\frac{1}{|f(z)-g(z)|}
$$

\end{proof}

\section{Approximation for complex curves in $\mathbb C^2$}

We assume that for every $c=(a,b)=(a+ib) \in \mathbb{C}$ we have a
holomorphic graph  $\Gamma_{c}$ given by $w=y_1+iy_2=f_{c}(z),
z=x_1+ix_2\in\Delta.$ We assume that all surfaces are disjoint and
that there is a surface through every point in
$\triangle\times\mathbb{C}.$ We assume that $f_{c}(0)=c.$

\medskip

Let $\pi:\Delta\times \mathbb{C}\rightarrow \mathbb{C}$ be defined
by $\pi(z,f_c(z))=c$.  The lamination of $\Delta \times \mathbb{C}$
by the $\Gamma_c$'s defines a holomorphic motion and so by
\cite{MSS} the map $(z,c)\mapsto(z,f_c(z))$  is a continuous
function in $(z,c)$.  It follows that the function $\pi$ is
continuous.

\medskip

Fix a positive constant $R$.  By Corollary \ref{estimate} there
exists a positive real number $\delta_0>0$  such that if
$z\in\frac{1}{2}\triangle$ and if $c,c'\in R\triangle$  with
$|c-c'|<\delta_0$ then

$$
(1) \left|\frac{\partial}{\partial z} f_{c'}(z)-
\frac{\partial}{\partial z}f_c(z)\right| \leq 4\cdot
|f_{c'}(z)-f_c(z)| \log \frac{1}{|f_{c'}(z)-f_c(z)|}.
$$

\medskip

We define a class of partially smooth functions:

\bea
{\mathcal A} & := & \{\phi \in {\mathcal C}(\Delta \times \mathbb{C}); \phi(z,f_c(z))\in {\mathcal C}^1(\Gamma_c),\\
& & \Phi(x_1,x_2,w):= \frac{\partial}{\partial x_1}\phi(x_1,x_2,f_c(x_1,x_2)), w=f_c(x_1,x_2) \in {\mathcal C}(\Delta \times \mathbb{C}),\\
& & \Psi(x_1,x_2,w):= \frac{\partial}{\partial
x_2}\phi(x_1,x_2,f_c(x_1,x_2)), w=f_c(x_1,x_2) \in {\mathcal
C}(\Delta \times \mathbb{C})\}. \eea

\begin{theorem}\label{approx}
Let $\phi\in \mathcal A$, let $R$ be a positive real number and let
$\epsilon >0$. Then there exists a function $\psi\in {\mathcal
C}^1(\triangle\times R\triangle)$ such that for every point
$(x_1,x_2,w)=(x_1,x_2,f_c(x_1,x_2))\in\triangle\times R\triangle$:
\bea
|\psi(x_1,x_2,w)-\phi(x_1,x_2,w)| & < & \epsilon,\\
|\frac{\partial}{\partial x_1}  [\psi(x_1,x_2,f_c(x_1,x_2))]-\frac{\partial}{\partial x_1} [\phi(x_1,x_2,f_c(x_1,x_2))]| & < & \epsilon,\\
|\frac{\partial}{\partial x_2}  [\psi(x_1,x_2,f_c(x_1,x_2))]-\frac{\partial}{\partial x_2} [\phi(x_1,x_2,f_c(x_1,x_2))]| & < & \epsilon.
\eea
\end{theorem}

We will prove the theorem using the following result:

\

\begin{proposition}\label{const}
Let $g\in \mathcal A, g(x_1,x_2,f_{a+ib}(x_1,x_2))=a$, and let $R$
be a positive real number.  There exists a positive real number
$t_0$ such that the following holds: For all $\epsilon
>0$ there exists a function $h\in {\mathcal
C}^1(t_0\triangle\times R\triangle)$ such that for every point
$(x_1,x_2,w)=(x_1,x_2,f_c(x_1,x_2))\in t_0\triangle\times
R\triangle$: \bea
|h(x_1,x_2,w)-g(x_1,x_2,w)| & < & \epsilon,\\
|\frac{\partial}{\partial x_1}  [h(x_1,x_2,f_c(x_1,x_2))]| & < & \epsilon,\\
|\frac{\partial}{\partial x_2}  [h(x_1,x_2,f_c(x_1,x_2))]| & < & \epsilon.
\eea
\end{proposition}

The same result holds if we replace $a$ by $b$ in the definition of $g$.

\medskip

\emph{Proof of Theorem \ref{approx} from Proposition \ref{const}:}
\

\begin{lemma}\label{impl}
Let $p\in\triangle$ be a point, and let $R,t_0$ be positive real
numbers such that $\triangle_{t_0}(p)\subset\subset\triangle$.
Consider the lamination restricted to
$\triangle_{t_0}(p)\times\mathbb{C}$. If the conclusion of
Proposition \ref{const} holds on $\triangle_{t_0}(p)\times
R\triangle$ (with respect to projection onto
$\{p\}\times\mathbb{C}$), then the conclusion of Theorem
\ref{approx} holds on $\triangle_{t_0}(p)\times R\triangle$.
\end{lemma}
\begin{proof}
Let $\pi=(\pi_1,\pi_2)$ denote the projection onto $\{p\}\times\mathbb{C}$.  For each $j,k\in\mathbb{Z}$ and $\delta>0$ we let $c^\delta(j,k)$ denote the point $(p,j\delta+k\delta i)$.  
Let $\Lambda^\delta_j$ denote the $\mathcal{C}^1$-smooth function defined by $\Lambda^\delta_j(t)=cos^2[\frac{\pi}{2\delta}(t-j\delta)]$ when $(j-1)\delta\leq t\leq (j+1)\delta$ and $0$ otherwise.  For each $c^\delta(j,k)$ we first define a function
$$
\psi^\delta_{jk}(z):=\phi(z,f_{c^\delta(j,k)}(z)),
$$  
and then we define a preliminary approximation
$$
\psi^\delta(z,w)=\sum_{j,k}\psi^\delta_{jk}(z)\Lambda_j(\pi_1(z,w))\Lambda_k(\pi_2(z,w)).
$$
Let $(z_0,w_0)\in\triangle_{t_0}(p)\times R\triangle$.  Then $\pi(z_0,w_0)$ is contained in a square with corners $c^\delta(j,k), c^\delta(j+1,k), c^\delta(j,k+1)$ and $c^\delta(j+1,k+1)$,  and we have that 
$$
\psi^\delta(z_0,w_0)=\sum_{m=j,j+1, n=k,k+1}\psi^\delta_{mn}(z_0)\Lambda_m(\pi_1(z_0,w_0))\Lambda_n(\pi_2(z_0,w_0)).
$$
We have that
\bea
& & |\psi^\delta(z_0,w_0)-\phi(z_0,w_0)|\\ 
& = & |\sum_{m=j,j+1, n=k,k+1} [\psi^\delta_{mn}(z_0)-\phi(z_0,w_0)]\cdot\Lambda^\delta_m(\pi_1(z_0,w_0)) \cdot\Lambda^\delta_n(\pi_2(z_0,w_0))|\\
& \leq & \mathrm{max}_{m=j,j+1, n=k,k+1}\{|\psi^\delta_{mn}(z_0)-\phi(z_0,w_0)|\}\\
\eea
Since the map from $\overline\triangle_{t_0}(p)\times\mathbb{C}$
defined by $(z,\alpha)\mapsto(z,f_\alpha(z))$ is a homeomorphism it follows that $\psi^\delta\rightarrow\phi$ uniformly as $\delta\rightarrow 0$.

\medskip

Next we approximate derivatives along leaves.  Let $\alpha$ be such that $(z_0,w_0)=(z_0,f_\alpha(z_0))$.  Since the functions $\Lambda^\delta_j\circ\pi_i$ are constant along leaves we get that
\bea
& & |\frac{\partial}{\partial x_i}[\psi^\delta(z_0,f_\alpha(z_0))-\phi(z_0,f_\alpha(z_0))]|\\ 
& = & |\sum_{m=j,j+1, n=k,k+1} [\frac{\partial}{\partial x_i}[\psi^\delta_{mn}(z_0)-\phi(z_0,f_\alpha(z_0))]]\\
& \times & \Lambda^\delta_m(\pi_1(z_0,f_\alpha(z_0))) \cdot\Lambda^\delta_n(\pi_2(z_0,f_\alpha(z_0)))|\\
& \leq & \mathrm{max}_{m=j,j+1, n=k,k+1}\{|\frac{\partial}{\partial x_i}[\psi^\delta_{mn}(z_0)-\phi(z_0,f_\alpha(z_0))]|\}\\
\eea
It follows that $\psi^\delta\rightarrow\phi$ also in $\mathcal{C}^1$-norm on leafs. 

\medskip

Now the conclusion of Lemma 2 follows because the functions $\pi_j$ can be approximated uniformly and in $\mathcal{C}^1$-norm on leaves. 
\end{proof}
For each point $p\in\triangle$  there exists by Proposition
\ref{const} a positive real number $t_p$ such that constant
approximation is possible on $\triangle_{t_p}(p)\times
R\triangle$.  Hence by Lemma \ref{impl} approximation of functions
in $\mathcal{A}$ is possible.

\medskip

We may then choose a locally finite cover
$\{U_\alpha\}_{\alpha\in\mathbb{N}}$ of $\triangle$ by disks such
that approximation by functions in $\mathcal{A}$ is possible on
each $U_\alpha\times R\triangle$.  Let $\{\varphi_\alpha\}$ be a
partition of unity subordinate to $\{U_\alpha\}$.  For each
$\alpha$ let $C_{\alpha}=\|\nabla\varphi_\alpha\|$.

\medskip

For a given $\epsilon_\alpha$ let $g_{\epsilon_\alpha}$ be an
$\epsilon_\alpha$-approximating function of $\phi$ on
$U_\alpha\times R\triangle$.  We will show that there is a
sequence $\{\epsilon_\alpha\}$ such that the function
$$
\psi=\sum_{\alpha}\varphi_\alpha\cdot\ g_{\epsilon_\alpha}
$$
satisfies the claims of the Theorem.

\medskip

Let $(z_0,f_c(z_0))\in U_\alpha$, and let $\alpha_1,...,\alpha_m$
be the finite number of $\alpha_i$'s such that the support of
$\phi_\alpha$ intersects $U_\alpha$.  Then
$$
\psi(z,f_c(z))=\sum_{i=1}^m\varphi_{\alpha_i}(z)\cdot\
g_{\epsilon_{\alpha_i}}(z,f_c(z))
$$
for all $z$ near $z_0$.  Then

\bea |\psi(z_0,f_c(z_0))-\phi(z_0,f_c(z_0))| & = & |[\sum_{i=1}^m
\varphi_{\alpha_i}(z_0)\cdot\
g_{\epsilon_{\alpha_i}}(z_0,f_c(z_0))]-\phi(z_0,f_c(z_0))|\\
& \leq & \sum_{i=1}^m\varphi_{\alpha_i}(z_0)\cdot|
g_{\epsilon_{\alpha_i}}(z_0,f_c(z_0))-\phi(z_0,f_c(z_0))|\\
& \leq & \mathrm{max}\{\epsilon_{\alpha_{i}}\}\\
\eea

Further

\bea
& & |\frac{\partial}{\partial x_1}[\psi(z_0,f_c(z_0))-\phi(z_0,f_c(z_0))]|\\
& = & |\frac{\partial}{\partial x_1}[[\sum_{i=1}^m
\varphi_{\alpha_i}(z)\cdot\
g_{\epsilon_{\alpha_i}}(z_0,f_c(z_0))]-\phi(z_0,f_c(z_0))]|\\
& = & |\sum_{i=1}^m \frac{\partial}{\partial
x_1}[\varphi_{\alpha_i}(z_0)\cdot(
g_{\epsilon_{\alpha_i}}(z_0,f_c(z_0))-\phi(z_0,f(z_0)))]|\\
& = & |\sum_{i=1}^m \frac{\partial}{\partial
x_1}[\varphi_{\alpha_i}(z_0)]\cdot(
g_{\epsilon_{\alpha_i}}(z_0,f_c(z_0)))-\phi(z_0,f(z_0)))\\
& + & \sum_{i=1}^m
\varphi_{\alpha_i}(z_0)\cdot\frac{\partial}{\partial x_1}[
g_{\epsilon_{\alpha_i}}(z_0,f_c(z_0))-\phi(z_0,f(z_0))]| \\
& \leq & m\cdot\mathrm{max}\{C_{\alpha_i}\}\cdot\mathrm{max}
\{\epsilon_{\alpha_i}\} + \mathrm{max}\{\epsilon_{\alpha_i}\}\\
\eea

Similarly we get that
$$
|\frac{\partial}{\partial
x_2}[\psi(z_0,f_c(z_0))-\phi(z,f_c(z_0))]|\leq
m\cdot\mathrm{max}\{C_{\alpha_i}\}\cdot\mathrm{max}
\{\epsilon_{\alpha_i}\} + \mathrm{max}\{\epsilon_{\alpha_i}\}
$$

It is clear that we may choose $\epsilon_{\alpha_i}$ for
$i=1,...,m$ to get the desired estimate for all points $z_0\in
U_\alpha$ for this particular $\alpha$.  Running through all
$\alpha$'s we have that any particular $\alpha_i$ will only come
under consideration a finite number of times.  Hence we may choose
the sequence $\{\epsilon_\alpha\}$. $\hfill\square$

\

We proceed to prove the Proposition.

Fix $\delta_0$ to get the estimate $(1)$ (in the beginning of Section 3) for all $|c-c'|<\delta_0$
with $|c|,|c'|\leq 2R$.  For any $\delta$  with $0<\delta<\delta_0$
we let $c^\delta(j,k)=(j+k\cdot i)\cdot\delta$ for $j,k \in \mathbb{Z}.$
Let $\chi: [0,1]\rightarrow \mathbb{R}$ be a smooth function such
that $\chi(t)=0$ for $0 \leq t \leq \frac{1}{4}$ and $\chi(t)=1$ for
$\frac{3}{4}\leq t \leq 1.$ Let $C$ be a constant such that
$|\chi'(t)| \leq C$ for all $t \in [0,1].$

\medskip

We first define a function $h_\delta$ on the surfaces
$\Gamma_{c^\delta(j,k)}$ simply by
${h_\delta}_{|\Gamma_{c^\delta(j,k)}}\equiv j\delta$. We want to
interpolate  this function between the surfaces.

\medskip

For a fixed $z$ consider the sets of points
$$
\Gamma_{c^\delta(j,k)}(z):=\{f_{c^{\delta}(j,k)}(z),f_{c^{\delta}(j+1,k)}(z),
f_{c^{\delta}(j,k+1)}(z),f_{c^{\delta}(j+1,k+1)}(z)\}.
$$
We first show that these sets move nicely with $z$  for small
enough $|z|$ and independent of $\delta$.  In particular we want
to know that we may define quadrilateral regions
$R_{\delta,j,k}(z)$ with corners $\Gamma_{c^\delta(j,k)}(z)$ and
that these sets have disjoint interior.

\medskip

We make the change of coordinates in the $w$ variable, by setting
$$
\tilde{w}(z,w)=\tilde{w}_{jk}(zw)=
\frac{w-f_{c^\delta(j,k)}(z)}{f_{c^\delta(j+1,k)}-f_{c^\delta(j,k)}(z)}.
$$

We get that

\bea
\tilde w(z,f_{c^\delta(j,k)}(z)) & \equiv & 0\\
\tilde w(z,f_{c^\delta(j+1,k)}(z)) & \equiv & 1\\
\eea

\begin{lemma}\label{grid}
Fix $N.$ Then there exists a real number $t_0>0$ independent of $\delta$  so that if
$|l|,|m|<N$ then $|\tilde w_{jk}(z,f_{c^\delta(j+l,k+m)}(z))- \tilde
w_{jk}(z,f_{c^\delta(j+l,k+m)}(0)) |<1/10 $ for all $|z|<t_0$ and any $j,k$.
\end{lemma}

\begin{proof}
Let $\pi:\Delta \rightarrow \mathbb C \setminus \{0,1\}$  be the
universal cover, and fix a pair $(j,k)$.  We have that
$p^{\delta}_{lm}:=\tilde w_{jk}(0,f_{c^\delta(j+l,k+m)}(0))=l+m\cdot i$  for
all $\delta$, so we may choose a point $P_{lm}\in\triangle$  such
that $\pi(P_{lm})=p^{\delta}_{lm}$ for all $\delta$.  For each
$\delta,j,k,l,m$ we have that $\tilde w_{jk}(z,f_{c^\delta(j+l.k+m)}(z))$  is a map
$k^{\delta}:\triangle\rightarrow\mathbb{C}\setminus\{0,1\}$, and so
they lift to maps $g^\delta:\triangle\rightarrow\triangle$  with
$g^\delta(0)=P_{lm}$, i.e. $k^\delta=\pi\circ g^\delta$.  By the
Schwartz' Lemma we have that $|g^\delta(z)-P_{lm}|\leq L_{lm}|z|$  for
all $\delta$.  So if $|z|$  is small enough we have that
$g^\delta(z)\subset\pi^{-1}(\triangle_{\frac{1}{10}}(p_{lm}))$.
Since there are only a finite number of pairs $(l,m)$ bounded by $N$
the result follows.
\end{proof}

From now on we assume that $|z|\leq t_0$.

\begin{lemma}
The quadrilaterals have disjoint interiors.
\end{lemma}

\begin{proof}
Pick $(j,k)$. We use the linear change of coordinates in the $w$
direction for fixed $z$:
$$
\tilde{w}_{jk}(z,w)=\frac{w-f_{c^{\delta}(j,k)}(z)}
{f_{c^{\delta}(j+1,k)}(z)-f_{c^{\delta}(j,k)}(z)}.
$$
This sends $f_{c^{\delta}(j+l,k+m)}(z)$ close to $(j+l,k+m)$ on
a small disc in the $z$ direction for uniformly bounded $(l,m)$. Hence
it is clear that the quadrilaterals are disjoint.
\end{proof}

Next we define preliminary functions $h^\delta_{jk}$ on the
respective quadrilaterals.  First we define a function
$t_z(y_1,y_2)$ to be constant equal to $0$ on the line between
$f_{c^\delta(j,k)}(z)$ and $f_{c^\delta(j,k+1)}(z)$, and constant
equal to $1$  on the line between $f_{c^\delta(j+1,k)}(z)$  and
$f_{c^\delta(j+1,k+1)}(z)$.  We extend $t_z$ continuously to be affine
on the two other edges, and then we extend $t_z$ to be constant equal
to $v$ on the line between $f_{c^\delta(j,k)}(z) +
v\cdot(f_{c^\delta(j+1,k)}(z)-f_{c^\delta(j,k)}(z))$ and
$f_{c^\delta(j,k+1)}(z) +
v\cdot(f_{c^\delta(j+1,k+1)}(z)-f_{c^\delta(j,k+1)}(z))$. Finally we
define $h^\delta_{jk}$  by
$$
h^\delta_{jk}(z,y_1,y_2)=j\delta + \delta\cdot(\chi\circ
t_z)(y_1,y_2).
$$
The $h^\delta_{jk}$'s patch up
smoothly in the "vertical" directions where the functions are
constant.  To be able to patch them together in the "horizontal"
directions we first extend each $h^\delta_{jk}$ across the
"horizontal" edges.

\medskip

To do this we use the coordinates defined by $\tilde w$.  Consider
the normalization
$$
\tilde
w_{jk}(z,w)=\frac{w-f_{c^\delta(j,k)}(z)}{f_{c^\delta(j+1,k)}(z)-f_{c^\delta(j,k)}(z)}
$$
Let $\tilde h_{jk}^\delta$ be defined by $\tilde h_{jk}^\delta\circ\tilde\omega = h_{jk}^\delta$. 
We want to glue together the two functions  on the quadrilaterals
sharing (in the new coordinates) the line segment $\gamma$ between
$(0,0)$ and $(1,0)$, i.e. the function $\tilde h^\delta_{jk}$ defined
above $\gamma$  and the function $\tilde h^\delta_{j(k-1)}$  below
$\gamma$.

\medskip

We start by extending the function $\tilde h^\delta_{jk}$.  Note first
that by Lemma \ref{grid} the quadrilaterals $R_{\delta,j,k}$  and
$R_{\delta,j,k-1}$ in the new coordinates - henceforth denoted $\tilde R_{\delta,j,k}$ and $\tilde R_{\delta,j,k-1}$ - have corners within
$\frac{1}{10}$-distance from the points $(l,m)$  for
$l,m\in\{0,1,-1\}$.  Note also that if we define a function
$\tilde t_z(\tilde y_1,\tilde y_2)$ ($\tilde w= \tilde y_1+i\tilde y_2$) along lines in the quadrilateral
$\tilde R_{\delta,j,k}(z)$ in the new coordinates as we did when we
defined $t_z(y_1,y_2)$  above, then
$h^\delta_{jk}=(j\delta+\delta(\chi\circ\tilde t))\circ\tilde w$.
Because of the placing of the corners we see that there exists a
constant $K$ independent of $\delta,j,k$  such that
$\|\nabla_{\tilde w}(j\delta+\delta(\chi\circ\tilde t))\|\leq K\delta$.

\medskip

Continue the lines in $\tilde R_{\delta,j,k}$ that pass through the
interval $[\frac{1}{8},1-\frac{1}{8}]$ and extend $\tilde h^\delta_{jk}$ to
be constant on these lines.  By the placing of the corners there is
a constant $\mu$ - independent of $\delta$ and $j,k$ - such that
these lines can be extended to the line between $(0,-\mu)$ and
$(1,-\mu)$.  Let $\tilde P_{\delta,j,k}$  denote the extended set
$\tilde R_{\delta,j,k}\cup(\tilde R_{\delta,j,k-1}\cap\{y_2\geq -\mu\})$, and we
see that $\tilde h^\delta_{jk}$  extends to be constant on the part of
$\tilde P_{\delta,j,k}$  where it is not already defined. Extend
$\tilde h^\delta_{j(k-1)}$  similarly in the other direction.

\medskip

To glue the functions together we choose a smooth function
$\varphi(z,\tilde y_1,\tilde y_2)=\varphi(\tilde y_2)$  such that $\varphi(\tilde y_2)=1$ if
$y_2\geq \mu$  and such that $\varphi(\tilde y_2)=0$  if $y_2\leq -\mu$.
We define
$$
h_\delta(z,w):=(\varphi\circ\tilde w_{jk})(z,w)\cdot
h^\delta_{jk}(z,w)+(1-\varphi\circ\tilde w_{jk})(z,w)\cdot
h^\delta_{j(k-1)}(z,w)
$$
Fix a constant $M$  such that $\|\frac{\partial\varphi}{\partial
\tilde y_2}\|= M$.

\begin{lemma}\label{change}
There are constants $N_1$ and $N_2$ such that for each $j,k,\delta$
we have that $h^\delta_{jk}(z,w)=j\delta$ if
$|w-f_{c^\delta(j,k)}(z)|\leq
N_1|f_{c^\delta(j+1,k)}(z)-f_{c^\delta(j,k)}(z)|$.  Moreover there
is a smooth function $\tilde g^\delta_{jk}(z,\tilde y_1,\tilde y_2)$ such that
$h^\delta_{jk}=\tilde g^\delta_{jk}\circ\tilde w$ and such that
$\|\nabla_{\tilde w}\tilde g^\delta_{jk}\|\leq N_2\delta$.
\end{lemma}
\begin{proof}
The existence of the constant $N_1$ can be seen by our description
of the function in local coordinates where we used Lemma \ref{grid}.
To see the rest let us give the function $\tilde g^\delta_{jk}$ explicitly.

\medskip

Fix $z$.  Let $(a_1,a_2)$  denote the corner of $\tilde R^\delta_{j,k}$  which
is close to $(0,1)$, and define a map
$A_z(\tilde y_1,\tilde y_2):=(\tilde y_1-\tilde y_2\frac{a_1}{a_2},\tilde y_2\frac{1}{a_2})$.  Then
$A_z$ changes smoothly with $z$ and we have that $\|A_z\|<2$ for all
the possibilities of $(a_1,a_2)$ we are considering.

\medskip

Next we define a function $\widehat t$ on the quadrilateral
$A_z(\tilde R^\delta_{j,k})$ along lines as above.  Let $(b_1,b_2)$  denote the
corner close to $(1,1)$ and fix $\widehat y=(\widehat y_1,\widehat y_2)$.  We have that the two
vertical sides of $A_z(\tilde R^\delta_{j,k})$ meet at the point $(0,-L)$ where
$L$ is given by $L=\frac{b_2}{b_1-1}$.  Calculating the slope of the
line from the point $\widehat y$ to the point $(\widehat t(\widehat y),0)$ we get that
$\frac{\widehat y_1}{L+\widehat y_2}=\frac{\tilde t(y)}{L}$ which  gives us
$$
\widehat t(\widehat y)=\frac{\widehat y_1\cdot L}{L+\widehat y_2}=\frac{\widehat y_1-b_2}{b_2+\widehat y_2(b_1-1)}.
$$
We have that $\widehat t$ varies smoothly with $(b_1,b_2)$   and we
see that $\widehat t$ has bounded derivatives for the cases of
$(b_1,b_2)$ we are considering.  Define $\tilde g^\delta_{jk}$ by
$$
\tilde g^\delta_{jk}=j\delta + \delta(\chi\circ\widehat t\circ A_z),
$$
and the function $h^\delta_{jk}$ is given by
$h^\delta_{jk}=\tilde g^\delta_{jk}\circ\tilde w$.
\end{proof}

\begin{lemma}
$h_\delta \rightarrow g$ in sup norm on $\triangle_{t_0}\times R\triangle$.
\end{lemma}

\begin{proof}
It is clear that $h_\delta(0,\cdot)\rightarrow g(0,\cdot)$
uniformly.  The claim then follows from Lemma  \ref{c1norm} below.
\end{proof}

\begin{lemma}\label{delta2}
If $t_0$ and $\delta$ is small enough we have that
$|f_{c^\delta(j,k)}(z)-f_{c^\delta((j+1,k)}(z)|\geq\delta^2$  for all
$z$  with $|z|\leq t_0$  and all $j,k$ such that
$|c^\delta(j,k)|\leq 2R$.
\end{lemma}

\begin{proof}
If $\delta$ is small enough we have the estimate
$$
|f_{c^\delta(j,k)}'(z)-f_{c^\delta(j+1,k)}'(z)|\leq 4|f_{c^\delta(j,k)}(z)-f_{c^\delta(j+1,k)}(z)|\log\frac{1}{|f_{c^\delta(j,k)}(z)-f_{c^\delta(j+1,k)}(z)|}
$$
for all $|z|\leq t_0$ and all $|c^\delta(j,k)|\leq 2R$ (Corollary \ref{estimate}).  Let $\phi$ denote the function 
$$
\phi(z)=|f_{c^\delta(j,k)}(z)-f_{c^\delta(j+1,k)}(z)|,
$$
and let $\phi(t)$ denote the restriction to radial real lines starting at the origin.  We have that
$$
\frac{\partial}{\partial t}[-\log\phi(t)]=-\frac{\phi'(t)}{\phi(t)},
$$
and so 
$$
|\frac{\partial}{\partial t}[-\log\phi(t)]|\leq4|-\log\phi(t)|.
$$
It follows that 
$$
-\log\phi(t)\leq-\log\phi(0)\cdot e^{4t}\Rightarrow\phi(t)\geq\phi(0)^{e^{4t}},
$$
and since $\phi(0)=\delta$  the result follows by choosing $t_0$  smaller than $\frac{1}{4}\log2$.
\end{proof}

\begin{lemma}\label{c1norm}
Let $c=a+ib, j\delta\leq a \leq (j+1)\delta, k\delta\leq b \leq
(k+1)\delta$. The function $h_\delta (z,f_c(z))$ is small in
${\mathcal C}^1$ norm along the graph $\Gamma_c$.
\end{lemma}

\begin{proof}

We need to estimate the derivatives of the function
$h_{\delta}(z,f_c(z))$ at an arbitrary point $(z_0,f_c(z_0))$.  We
estimate $\frac{\partial}{\partial x}= \frac{\partial}{\partial x_1}$ - the case of
$\frac{\partial}{\partial x_2}$ is similar.  Since we are working on
lines we use the notation $(x,y_1,y_2)$ for coordinates.

\medskip

We observe first that if $(z_0,f_c(z_0))$ is outside $R_{\delta,j,k}$ then it must
still be very close. If the point is close to the vertical edges, then the function $h_\delta$
is locally constant, so we are done. We can assume that also $(z_0,f_c(z_0))\in P_{\delta,j,k}
\setminus P_{\delta,j,k+1}$. We divide the proof
into two cases: Assume first that $(z_0,f_c(z_0))$ is not in
 $R_{\delta,j,k-1}.$  Then the function $h_\delta$  is simply
equal to the function $h^\delta_{jk}$.

\medskip

We have that

\bea \frac{\partial}{\partial x}(h^{\delta}_{jk}(x,f(x))) & = &
(\frac{\partial h^{\delta}_{jk}}{\partial x},\frac{\partial
h^{\delta}_{jk}}{\partial y_1},\frac{\partial
h^{\delta}_{jk}}{\partial y_2})(x,f(x))\cdot(1,\frac{\partial
f_1}{\partial x},\frac{\partial f_2}{\partial x})(x)\\
& = & \frac{\partial h^{\delta}_{jk}}{\partial x}(x,f(x)) +
(\frac{\partial h^{\delta}_{jk}}{\partial y_1},\frac{\partial
h^{\delta}_{jk}}{\partial y_2})(x,f(x))\cdot(\frac{\partial
f_1}{\partial
x},\frac{\partial f_2}{\partial x})(x)\\
\eea

For fixed $s,v$ we may define a curve $(x,g(x))$:

\bea g(x) & = & (1-s)[(1-v)f_{c^\delta(j,k)}(x) +
vf_{c^\delta(j+1,k)}(x)]\\
& + & s[(1-v)f_{c^\delta(j,k+1)}(x) +
vf_{c^\delta(j+1,k+1)}(x)].\\
\eea

Then $h^{\delta}_{jk}(x,g(x))\equiv j \delta+\chi(v) \delta$.  Choose $s$
and $v$ so that $(x_0,g(x_0))=(x_0,f_c(x_0))$. We get that

$$
0=\frac{\partial}{\partial x}(h^{\delta}_{jk}(x,g(x)))=\frac{\partial
h^{\delta}_{jk}}{\partial x}(x,g(x)) + (\frac{\partial
h^{\delta}_{jk}}{\partial y_1},\frac{\partial
h^{\delta}_{jk}}{\partial y_2})(x,g(x))\cdot(\frac{\partial
g_1}{\partial x},\frac{\partial g_2}{\partial x})(x),
$$
and so
$$
\frac{\partial}{\partial
x}(h^{\delta}_{jk}(x_0,f(x_0)))=(\frac{\partial
h^{\delta}_{jk}}{\partial y_1},\frac{\partial
h^{\delta}_{jk}}{\partial y_2})(x_0,g(x_0))\cdot(\frac{\partial
f_1}{\partial x}-\frac{\partial g_1}{\partial x},\frac{\partial
f_2}{\partial x}-\frac{\partial g_2}{\partial x})(x_0).
$$
Using the Lemma \ref{grid} we see that
$\|f_c(x_0)-f_{c^\delta(j+l,k+m)}(x_0)\|\leq
2\|f_{c^\delta(j+1,k)}(x_0)-f_{c^\delta(j,k)}(x_0)\|$ for
$l,m\in\{0,1\}$, and so

\bea & & \|\frac{\partial}{\partial
x}(f_c-f_{c^\delta(j+l,k+m)})(x_0)\|\\
& \leq & 4\|(f_c-f_{c^\delta(j+l,(k+m)})(x_0)\|\log\frac{1}{\|
(f_c-f_{c^\delta(j+l,k+m)})(x_0)\|}\\
& \leq & \
8\|(f_{c^\delta(j+1,k)}-f_{c^\delta(j,k)})(x_0)\|\log\frac{1}{2\|
(f_{c^\delta(j+1,k)}-f_{c^\delta(j,k)})(x_0)\|}.\\
\eea

It follows that

\bea & & \|\frac{\partial}{\partial x}(h_{\delta}(x,f(x))\|\\
& \leq & 8\cdot\|(\frac{\partial h_{\delta}}{\partial
y_1},\frac{\partial h_\delta}{\partial
y_2})\|\cdot\|(f_{c^\delta(j+1,k)}-f_{c^\delta(j,k)})(x_0)\|\\
& \times & \log\frac{1}{2\|
(f_{c^\delta(j+1,k)}-f_{c^\delta(j,k)})(x_0)\|}. \eea

We proceed to estimate $\|(\frac{\partial h_{\delta}}{\partial
y_1},\frac{\partial h_{\delta}}{\partial y_2})\|$.  We change
coordinates according to Lemma \ref{change} and write $h_{\delta}$
as a composition $\tilde g_\delta\circ\tilde w(y)$.  We get $\|D_w\tilde
w\|=\frac{1}{\|f_{c^\delta(j+1,k)}(x_0)-f_{c^\delta(j,k)}(x_0)\|}$,
and we have that $\|\nabla_{\tilde w}\tilde g_\delta\|\leq N_2\delta$. This
shows that

$$
\|(\frac{\partial h_{\delta}}{\partial y_1},\frac{\partial
h_{\delta}}{\partial y_2})\|\leq
N_2\delta\frac{1}{\|f_{c^\delta(j+1,k)}(x_0)-f_{c^\delta(j,k)}(x_0)\|}.
$$

This gives
$$
\|\frac{\partial}{\partial x}(h_{\delta}(x,f(x))\|\leq
8N_2\delta\log\frac{1}{\|f_{c^\delta(j+1,k)}(x_0)-f_{c^\delta(j,k)}(x_0)\|}.
$$
We have by Lemma \ref{delta2} that
$\|f_{c^\delta(j+1,k)}(x_0)-f_{c^\delta(j,k)}(x_0)\|\geq \delta^2$
and so
$$
\|\frac{\partial}{\partial x}(h_{\delta}(x_0,f(x_0))\|\leq
8N_2\delta\log\frac{1}{2\delta^2}\rightarrow 0 \ \mathrm{as} \
\delta\rightarrow 0.
$$

\medskip

The other case we have to consider is when $(z_0,f_c(z_0))$  is
contained in an overlap where we glued our functions together.  In
that case we may assume that $(z_0,f_c(z_0))$ is also contained in
$P^\delta_{j(k-1)}$.

\

Let $\overrightarrow{v}$  denote the vector
$\overrightarrow{v}=\frac{\partial}{\partial x}(x_0,f_c(x_0))$. We
have that

\

\bea
\nabla h_\delta(x_0,f_c(x_0))\cdot\overrightarrow{v} & = &
\nabla[\varphi\circ\tilde w\cdot
h^\delta_{jk}](x_0,f_c(x_0))\cdot\overrightarrow{v}\\
& + & \nabla[(1-\varphi)\circ\tilde w\cdot
h^\delta_{j(k-1)}](x_0,f_c(x_0))\cdot\overrightarrow{v}\\
& = & h^\delta_{jk}(x_0,f_c(x_0))\cdot \nabla[\varphi\circ\tilde
w](x,f_c(x_0))\cdot\overrightarrow{v}\\
& + & (\varphi\circ\tilde w)(x_0,f_c(x_0))\cdot
\nabla[h^\delta_{jk}](x_0,f_c(x_0))\cdot\overrightarrow{v}\\
& + & h^\delta_{j(k-1)}(x_0,f_c(x_0))\cdot
\nabla[(1-\varphi)\circ\tilde
w](x,f_c(x_0))\cdot\overrightarrow{v}\\
& + & ((1-\varphi)\circ\tilde w)(x_0,f_c(x_0))\cdot
\nabla[h^\delta_{j(k-1)}](x_0,f_c(x_0))\cdot\overrightarrow{v}\\
\eea

By the above calculations we need not worry about the second and
fourth term in this sum so we have to check that

$$
(h^\delta_{jk}(x_0,f_c(x_0))-h^\delta_{j(k-1)}(x_0,f_c(x_0)))\cdot
\nabla[\varphi\circ\tilde
w](x_0,f_c(x_0))\cdot\overrightarrow{v}\rightarrow 0
$$

\

as $\delta\rightarrow 0$.

\medskip

First of all we have that
$|h^\delta_{jk}(x_0,f_c(x_0))-h^\delta_{j(k-1)}(x_0,f_c(x_0))|\leq
2\delta$. Further $|\nabla[\varphi\circ\tilde
w](x_0,f_c(x_0))\cdot\overrightarrow{v}|\leq M\cdot\|D[\tilde
w](x_0,f_c(x_0))(\overrightarrow{v})\|$. \

Now
$$
D[\tilde
w](x_0,f_c(x_0))(\overrightarrow{v})=\frac{\partial}{\partial
x}[(x,\frac{f_c(x)-f_{c^\delta(j,k)}(x)}{f_{c^\delta(j+1,k)}(x)
-f_{c^\delta(j,k)}(x)})](x_0).
$$
Ignoring the constant term (it gets killed by $\delta$) we get that

\

\bea & & \|D[\tilde w](x_0,f_c(x_0))(\overrightarrow{v})\|\\ &
\leq & \frac{|f_c'(x_0)-f_{c^\delta(j,k)}'(x_0)|}{|f_{c^\delta(j+1,k)}(x_0)-f_{c^\delta(j,k)}(x_0)|}\\
& + &
\frac{|f_{c}(x_0)-f_{c^\delta(j,k)}(x_0)|\cdot|f_{c^\delta(j+1,k)}'(x_0)-f_{c^\delta(j,k)}'(x_0))|}
{|f_{c^\delta(j+1,k)}(x_0)-f_{c^\delta(j,k)}(x_0)|^2}\\
& \leq  & \frac{|f_c(x_0)-f_{c^\delta(j,k)}(x_0)|}{|f_{c^\delta(j+1,k)}(x_0)-f_{c^\delta(j,k)}(x_0)|}\log\frac{1}{|f_c(x_0)-f_{c^\delta(j,k)}(x_0)|}\\
& + &
\frac{|(f_{c}(x_0)-f_{c^\delta(j,k)}(x_0))|\cdot|(f_{c^\delta(j+1,k)}(x_0)-f_{c^\delta(j,k)}(x_0))|}{|(f_{c^\delta(j+1,k)}(x_0)-f_{c^\delta(j,k)}(x_0))|^2}\\
& \times & \log\frac{1}{|(f_{c^\delta(j+1,k)}(x_0)-f_{c^\delta(j,k)}(x_0))|}\\
\eea

\medskip

By Lemma \ref{grid} we have that
$\frac{|f_c(x_0)-f_{c^\delta(j,k)}(x_0)|}{|f_{c^\delta(j+1,k)}(x_0)-f_{c^\delta(j,k)}(x_0)|}\leq
2$ and so

\bea \|D[\tilde w](x_0,f_c(x_0))(\overrightarrow{v})\| & \leq &
2\cdot\log\frac{1}{|f_{c^\delta(j+1,k)}(x_0)-f_{c^\delta(j,k)}(x_0)|}\\
& + & 2\log\frac{1}{|f_c(x_0)-f_{c^\delta(j,k)}(x_0)|}\\
\eea

\

By Lemma \ref{change} we have that our function is constant unless
$|f_c(x_0)-f_{c^\delta(j,k)}(x_0)|\geq
N_1|f_{c^\delta(j+1,k)}(x_0)-f_{c^\delta(j,k)}(x_0)|\geq
N_1\delta^2$ (by Lemma 7), and so we may assume that

$$
\|D[\tilde w](x_0,f_c(x_0))(\overrightarrow{v})\|\leq
2\log\frac{1}{\delta^2}+2\log\frac{1}{N_1\delta^2}.
$$

\

All in all:

\bea & &
|(h^\delta_{jk}(x_0,f_c(x_0))-h^\delta_{j(k-1)}(x_0,f_c(x_0)))\cdot
\nabla[\varphi\circ\tilde
w](x_0,f_c(x_0))\cdot\overrightarrow{v}|\\
& \leq &  4M\delta
(\log\frac{1}{\delta^2}+\log\frac{1}{N_1\delta^2})\rightarrow 0 \ \mathrm{as}\ \delta\rightarrow 0.\\
\eea

\

\end{proof}

\section{Proof of the main theorem}

We are ready to prove the main theorem.
By the theorem of Slodkowski, \cite{Sl1991}, \cite{AM2001}, we can assume that $\mathcal{L}$ is a lamination of $\triangle\times\mathbb{C}$  as in the previous section. 

\medskip

\emph{Proof of the Main Theorem:}
Suppose that $T$ is a positive closed $(1,1)$ current on $\Delta^2(0,1)$, supported on the laminated set $K$ described in the introduction. We assume that $T$ is directed by the lamination $\mathcal L$ of $K.$ Hence there is a positive measure $\mu$ so that $T=\int [V_\alpha]d\mu(\alpha).$ Suppose that
 $\lambda=dw-f_\alpha'(z)dz$. We want to show that $\lambda \wedge T=0.$
 Let $\phi$ be any smooth $(1,0)$ test form. We need to show that $<\lambda\wedge T,\phi>=0.$
 This follows since 
 
 \bea
 <\lambda \wedge T,\phi> & = & \int(\lambda \wedge T) \wedge \phi\\
 & = & \int T \wedge (\lambda \wedge \phi)\\
 & = & \int_\alpha \left(\int_{V_\alpha} \lambda \wedge \phi\right) d\mu(\alpha)\\
 & = & \int_\alpha 0=0\\
 \eea
 
\medskip

Assume next that $T$ is weakly directed by $\mathcal{L}$.  Since $\mathcal{L}$ is a lamination of $\triangle\times\mathbb{C}$ we may invoke the approximation result from the previous section.  With the approximation result at hand the implication follows from 
Sullivan's proof of the smooth case \cite{S1976}.  We include the proof for the benefit of the reader. 

\medskip

Step 1 is to show that there exists a family of probability measures $\sigma_\alpha$ such that $\sigma_\alpha$ is supported on $\Gamma_\alpha$, and a measure $\mu'$ on the $\alpha$-plane, such that for all test forms $\omega$ we have that 
$$
T(\omega)=\int ( \int_{\Gamma_\alpha}\omega d\sigma_\alpha ) \ d\mu'.
$$
Let $\omega$ be a $(1,1)$ test form and let $\lambda(z,w)=dw-f'_\alpha(z)dz$ for $w=f_\alpha(z)$.  Let $\overrightarrow{v_1}(z,w)=(1,f'_\alpha(z))$  and let  $\overrightarrow{v_2}(z,w)=(i,i\cdot f'_\alpha(z))$
for $w=f_\alpha(z)$, and define the 2-tangent field $v(z,w)=(\overrightarrow{v_1}(z,w),\overrightarrow{v_2}(z,w))$.  

\medskip

Switching basis we have that 
$$
\omega=\psi_1 dz\wedge d\overline z + \psi_2 dz\wedge\overline\lambda + \psi_3 d\overline z\wedge\lambda + \psi_4 \lambda\wedge\overline\lambda
$$
for some functions $\psi_i$, and by assumption we have that $T(\omega)=T(\psi_1dz\wedge d\overline z)$.  The function $\psi_1$ is given by $\psi_1=\frac{1}{2i}\omega(v)$ and so we have that 
$$
T(\omega)=T(\frac{1}{2i}\omega(v)dz\wedge d\overline z).
$$
On the other hand we may use $T$ to define a linear functional $L$ on $\mathcal{C}_0(\triangle\times\mathbb{C})$ by  $L(\psi)=T(\psi dz\wedge d\overline z)$, and so by Riesz' Representation Theorem there is a measure $\nu$ such that 
$$
L(\psi)=\int\psi d\nu.
$$
This means that 
$$
T(\omega)=\int\frac{1}{2i}\omega(v)d\nu.
$$

\medskip

Now the measure $\nu$ disintegrates \cite{HA}: There exists a family of probability measures $\sigma_\alpha$ such that $\sigma_\alpha$ is supported on $\Gamma_\alpha$, and  a measure $\mu'$ on the $\alpha$-plane, such that for all $\psi\in\mathcal{C}_0(\triangle\times\mathbb{C})$  we have that 
$$
\int\psi d\nu = \int  ( \int_{\Gamma_\alpha} \psi d\sigma_\alpha ) \ d\mu'.
$$

We define currents $T_\alpha$  by $T_\alpha(\omega)=\int_{\Gamma_\alpha}\frac{1}{2i}\omega(v)d\sigma_\alpha$, and we get that 
$$
T(\omega)=\int T_\alpha(\omega)d\mu'.
$$

\medskip

The next step is to show that $T_\alpha$ is closed for $\mu'$-almost all $\alpha$.  Let $\{\omega_j\}$ be a dense set of $\mathcal{C}^1$-smooth $(0,1)$ test forms and fix a $j\in\mathbb{N}$.  Let $g$ be a continuous function in the $\alpha$-variable and extend $g$ constantly along leaves.  We want to show that 
$$
\int g\cdot T_\alpha(\partial\omega)d\mu' = 0
$$
because this would imply that $\partial T_\alpha=0$ for $\mu'$-almost all $\alpha$ (since $g$ is arbitrary).

\medskip

By Theorem 1 there exists a sequence $g_i$ of smooth functions such that $g_i\rightarrow g$  uniformly and in $\mathcal{C}^1$-norm on leaves.  Since $T$ is closed we have that 
$$
0=\int T_\alpha(\partial(g\omega_j))d\mu'=\int T_\alpha(\partial g_i\wedge\omega_j)d\mu' + \int g_i\cdot T_\alpha(\partial\omega_j)d\mu'. 
$$
Since $T_\alpha(\partial g_i\wedge\omega)\rightarrow 0$  we get that 
$$
\int g\cdot T_\alpha(\partial\omega_j)d\mu'=\lim_{i\rightarrow\infty}\int g_i\cdot T_\alpha(\partial\omega_j)d\mu' = 0 
$$ 
Running through all $\omega_j$'s we see that $T_\alpha$ is closed for $\mu'$-almost all $\alpha$.  The only possibility then is that the measures $\sigma_\alpha$  are constant multiples of $dz\wedge d\overline z$, i.e. $\sigma_\alpha=\varphi(\alpha)dz\wedge d\overline z$  where $\varphi$ is a measurable function \cite{LE}.  Define $\mu:=\varphi\cdot \mu'$.

\section{Two Counterexamples}

In \cite{FWW2007} the authors proved versions of the main theorem for real laminations
in $\mathbb R^2$ and $\mathbb R^3.$ In those results we added an extra slope condition on the
laminations which is analogous to the estimate in Corollary 1. We give here a simple example
of a lamination of curves in $\mathbb R^2$ where the slope condition is not satisfied.
Also the conclusion of the Main Theorem fails. The analogue of Theorem 1, i.e. approximation
of partially smooth functions fails as well.

\bigskip

For each $t\in \mathbb R$, we let $\gamma_t$ be the curve $y=f_t(x)=(x-t)^3$ in $\mathbb R^2.$
Clearly this gives a continuous lamination of $\mathbb R^2$ by curves. The curves are all tangent to the $x-$ axis. This implies that the current of integration of the $x$- axis is annihilated by the one form
$\lambda$ which is defined by $dy-f'_t(x)dx$ on $\gamma_t.$ However, this current is not an integral of currents $[\gamma_t].$ We also observe that the function $a(x,y)$ defined by
$a(x,f_t(x))=t$ cannot be approximated by ${\mathcal C}^1$ functions, because any such approximation 
will have to have a small derivative along the $x-$ axis.

\bigskip

We can also modify this example so that we have a Riemann surface lamination in $\mathbb C^3.$
For $t\in \mathbb C$, let $\gamma_t$ be the complex curve $\gamma_t(s)=(z,w,\tau)=(s,(s-t)^2,(s-t)^3).$
These curves laminate $\mathbb C^3$ and $\gamma_t$ it tangent to the $z-$ axis at $(t,0,0).$
Hence the $z-$ axis is annihilated by any continuous one forms defining the lamination. Hence the
current of integration of the $z$ axis is weakly directed. But clearly it is not directed by the lamination.
Again the function $a(z,w,\tau)$ defined by $a_{|\gamma_t}=t$ cannot be approximated by $\mathcal C^1$ functions.

\bigskip

\noindent John Erik Forn\ae ss\\
Mathematics Department\\
The University of Michigan\\
East Hall, Ann Arbor, MI 48109\\
USA\\
fornaess@umich.edu\\

\noindent Yinxia Wang\\
Department of Mathematics\\
Henan University\\
Kaifeng, 475001\\
China\\
yinxiawang@gmail.com\\

\noindent Erlend Forn\ae ss Wold\\
Mathematisches Institut\\
Universit\"at Bern\\
Sidlerstr. 5\\
CH-3012 Bern\\
Switzerland\\
erlendfw@math.uio.no\\

\end{document}